\documentclass[aos,MSNbibl,nameyear,seceqn,dvips]{arximspdf}

\doi{10.1214/14-AOS1212} 
\volume{42}
\issue{3}
\pubyear{2014}
\firstpage{1003}
\lastpage{1028}

\makeatletter
\newcommand{\rrvert}{\vert}
\newcommand{\llvert}{\vert}
\def\implies{\Rightarrow}
\newtheorem{theorem}{Theorem}
\newtheorem{corollary}{Corollary}
\newproclaim{definition}{Definition}
\newproclaim{example}{Example}
\newtheorem{lemma}{Lemma}
\newtheorem{proposition}{Proposition}
\makeatother

\begin{document}
\begin{frontmatter}

\title{Merging and testing opinions}
\runtitle{Merging and testing opinions}

\begin{aug}
\author{\fnms{Luciano} \snm{Pomatto}\corref{}\ead[label=e2]{l-pomatto@kellogg.northwestern.edu}\ead[label=u2,url]{http://www.kellogg.northwestern.edu/faculty/pomatto/index.htm}},
\author{\fnms{Nabil} \snm{Al-Najjar}\ead[label=e1]{al-najjar@kellogg.northwestern.edu}\ead[label=u1,url]{http://www.kellogg.northwestern.edu/faculty/directory/al-najjar\_nabil.aspx}}
\and
\author{\fnms{Alvaro} \snm{Sandroni}\ead[label=e3]{sandroni@kellogg.northwestern.edu}\ead[label=u3,url]{http://www.kellogg.northwestern.edu/faculty/directory/sandroni\_alvaro.aspx}\thanksref{t1}}

\runauthor{L. Pomatto, N. Al-Najjar and A. Sandroni}
\thankstext{t1}{Supported by a grant from the NSF.}
\affiliation{Northwestern University}

\address{Managerial Economics \& Decision Sciences Department\\
Kellogg School of Management\\
2001 Sheridan Road\\
Evanson, Illinois 60208\\
USA\\
\printead{e1}\\
\phantom{E-mail:\ }\printead*{e2}\\
\phantom{E-mail:\ }\printead*{e3}\\
\printead{u1}\\
\phantom{URL:\ }\printead*{u2}\\
\phantom{URL:\ }\printead*{u3}}
\end{aug}

\received{\smonth{8} \syear{2013}}
\revised{\smonth{12} \syear{2013}}

%
\begin{abstract}
We study the merging and the testing of opinions in the context of a
prediction model. In the absence of incentive problems, opinions can be
tested and rejected, regardless of whether or not data produces consensus
among Bayesian agents. In contrast, in the presence of incentive problems,
opinions can only be tested and rejected when data produces consensus among
Bayesian agents. These results show a strong connection between the testing
and the merging of opinions. They also relate the literature on Bayesian
learning and the literature on testing strategic experts.
\end{abstract}

%


\end{frontmatter}

\section{Introduction}\label{sec1}

Data can produce consensus among Bayesian agents who initially
disagree. It
can also test and reject opinions. We relate these two critical uses of data
in a model where agents may strategically misrepresent what they know.

In each period, either $0$ or $1$ is observed. Let $P$ and $Q$ be two
probability measures on $ \{ 0,1 \} ^{\infty}$ such that $Q$ is
absolutely continuous with respect to $P$. If $P$ and $Q$ are $\sigma$-additive
then, as shown by \citet{BlaDub62}, the conditional
probabilities of $P$ and $Q$ merge, in the sense that the two posteriors
become uniformly close as the amount of observations increases ($Q$-almost
surely). So, repeated applications of Bayes' rule lead to consensus among
Bayesian agents, provided that their opinions were initially compatible.

Now consider Savage's axiomatization of subjective probability. He proposed
postulates that characterize a preference relation over bets
in terms of a nonatomic finitely additive probability $P$. Call such $P$,
for short, an \textit{opinion}. Savage's framework allows for finitely
additive probabilities that are not $\sigma$-additive. In particular,
the conclusions of the Blackwell and Dubins theorem hold for some, but
not all, opinions. This
flexibility makes Savage's framework an ideal candidate to study the
connection between the merging and the testing of opinions.

We say that an opinion $P$ satisfies the \textit{Blackwell--Dubins property}
if whenever $Q$ is an opinion absolutely continuous with respect to
$P$, the
two conditional probabilities merge. By definition, in this subframework,
sufficient data produces agreement among Bayesian agents who have compatible
initial opinions. Outside this subframework, Bayesian agents may satisfy
Savage's axioms, have compatible initial opinions and yet persistently
disagree. See the \hyperref[app]{Appendix} for an example.

Any opinion, whether or not it satisfies the Blackwell and Dubins property,
can be tested and rejected. To reject an opinion $P$, it suffices to
find an
event that has low probability according to $P$ and then reject it if this
event is observed. Thus, if opinions are honestly reported then the
connection between merging and testing opinions is weak. In the absence of
incentive problems, subjective probabilities can be tested and rejected
whether or not data produces consensus.

Now consider the case in which a self-proclaimed expert, named Bob, may
strategically misrepresent what he knows. Let Alice be a tester who
wants to
determine whether Bob is an informed expert who honestly reports what he
believes or he is an uninformed, but strategic, expert who has reputational
concerns and wants to pass Alice's test. Alice faces an adverse selection
problem and uses data to screen the two types of experts.

A test is likely to \textit{control for type I error} if an informed expert
expects to pass the test by truthfully reporting what he believes. A test
can be \textit{manipulated} if even completely uninformed experts are likely
to pass the test, no matter how the data unfolds in the future. The word
``likely'' refers to a possible
randomization by the strategic expert to manipulate the test. Only
nonmanipulable tests that control for type I error pass informed experts
and may fail uninformed ones.

Our main results are: In the presence of incentive problems, if opinions
must satisfy the Blackwell--Dubins property then there exists a test that
controls for type~I error and cannot be manipulated. If, instead, any
opinion is allowed then every test that controls for type I error can be
manipulated. Thus, in Savage's framework strategic experts cannot be
discredited. However, strategic experts can be discredited if opinions are
restricted to a subframework where data produces consensus among Bayesian
agents with initially compatible views. These results show a strong
connection between the merging and the testing of opinions but only under
incentive problems.

The Blackwell--Dubins theorem has an additional interpretation. In this
interpretation, $Q$ is referred to as the data generating process and
$P$ is
an agent's belief initially compatible with $Q$. When the conclusions
of the
Blackwell--Dubins theorem hold, then $P$ and $Q$ merge and so, the agent's
predictions are eventually accurate. Thus, multiple repetitions of Bayes'
rule transforms the available evidence into a near perfect guide to the
future. It follows that our main results also have an additional
interpretation. Under incentive problems, strategic experts can only be
discredited if they are restricted to a subframework where opinions
that are
compatible with the data generating process are eventually accurate.

Finally, our results relate the literatures on Bayesian learning and the
literature on testing strategic experts (see the next section for
references). They show a strong connection between the framework under which
Bayesian learning leads to accurate opinions and the framework under which
strategic experts can be discredited.

The paper is organized as follows. Section~\ref{sec2} describes the model.
Section~\ref{sec3}
reviews the Blackwell--Dubins theorem and defines the Blackwell--Dubins
property. Section~\ref{sec4} contains our main results. Section~\ref{sec5} relates our results
and category tests. Section~\ref{sec6} considers the case where the set of
per-period outcome may be infinite. The \hyperref[app]{Appendix} contains all proofs
and a
formal example of a probability that does not satisfy the Blackwell--Dubins
property.

\subsection{Related literature}\label{sec1.1}

Blackwell and Dubins' idea of merging of opinions is central in the
theory of Bayesian learning
and Bayesian statistics. In Bayesian nonparametric statistics, see the
seminal work of \citet{DiaFre86}, \citet{DArDiaFre88} and the more
recent work by Walker, Lijoi and
Pruenster (\citeyear{WalLijPru05}). In the theory of Bayesian learning, see \citet{SchSei90}. We refer to \citet{Daw85} for a connection with the
theory of calibration.

In game theory, the Blackwell--Dubins theorem is central in the study
of convergence to Nash equilibrium in repeated games. The main
objective is to understand the conditions under which
Bayesian learning leads to a Nash equilibrium [see, among many
contributions, Foster and Young (\citeyear{FosYou01}, \citeyear{FosYou03}),
\citet{FudKre93},
Fudenberg and Levine (\citeyear{FudLev98,FuLe09}), \citet{HarMas13},
\citet{JacKalSmo99}, Kalai and Lehrer (\citeyear{KalLeh93N1,KalLeh93N2}),
Lehrer and Smorodinsky (\citeyear{LehSmo96N1,LehSmo96N2}), \citet{MonSamSel97},
Nachbar (\citeyear{Nac97}, \citeyear{Nac01},
\citeyear{Nac05}), \citet{San98} and Young (\citeyear{You02}, \citeyear{You04})].

A series of papers investigate whether empirical tests can be
manipulated. In statistics, see \citet{FosVoh98}, Cesa-Bianchi
and Lugosi (\citeyear{CesLug06}), \citet{VovSha05} and
Olszewski and Sandroni
(\citeyear{OlsSan09N1}). In economics, see among several contributions,
\citet{AlNWei08},
Al-Najjar et al. (\citeyear{AlNetal10}), Babaioff et al. (\citeyear{Babetal11}), \citet{DekFei06}, \citet{FeiLam}, \citet{FeiSte08},
\citet{ForVoh09}, \citet{FudLev99}, \citet{GraSal},
\citet{GraShm}, Hu and Shmaya (\citeyear{HuShm13}), \citet{Leh01},
Olszewski and Peski (\citeyear{OlsN1}), Olszewski and
Sandroni (\citeyear{OlsSan07,OlsSan08,OlsSan09N1,OlsSan09N2,OlsSan11}), \citet{San03},
\citet{SanSmoVoh03},
Shmaya (\citeyear{Sh08}),\vadjust{\goodbreak}
\citet{Ste11}.
For a review, see \citet{FosVoh} and
\citet{OlsN2}. See also \citet{AlNPomSan} for a
companion paper.

\section{Setup}\label{sec2}

In every period an outcome, $0$ or $1$, is observed (all results generalize
to the case of finitely many outcomes). A \emph{path} is an infinite
sequence of outcomes and $\Omega= \{ 0,1 \} ^{\infty}$ is
the set
of all paths. Given a path $\omega$ and a period $t$, let $\omega
^{t}\subseteq\Omega$ be the cylinder of length $t$ with base $\omega$.
That is, $\omega^{t}$ is the set of all paths which coincide with
$\omega$
in the first $t$ periods. The set of all paths $\Omega$ is endowed
with a $%
\sigma$-algebra of \emph{events} $\Sigma$ containing all cylinders.

The set $\Omega$ is endowed with the product topology. In this
topology, a
set is open if and only if it is a countable union of cylinders. We
denote by $\Sigma_{1}$ the set of all open subsets of $\Omega$
and by $\mathcal{B}$ the Borel $\sigma$-algebra generated by the topology. Note that
$\Sigma_{1}\subset
\mathcal{B
}\subseteq\Sigma$.

Let $\mathbb{P}$ be the set of all finitely additive probabilities on $%
( \Omega,\Sigma ) $. A probability $P\in\mathbb{P}$ is strongly
nonatomic, or \textit{Savagean}, if for every event $E$ and every
$\alpha
\in [ 0,1 ] $ there is an event $F\subseteq E$ such that
$P (
F ) =\alpha P ( E ) $. The term ``Savagean'' emphasizes the relation between strongly
nonatomic probabilities and the \citet{Sav54} representation theorem: a finitely
additive probability corresponds to a preference relation satisfying Savage's
axioms if and only if it is strongly nonatomic. To simplify the
language, we also refer to a Savagean probability as an \textit{opinion}.
Let $\Delta$ denotes the set of all opinions.

At time $0$, a self-proclaimed expert, named Bob, announces an opinion $P$.
A~tester, named Alice, evaluates his opinion empirically. Alice announces
her test at period $0$, before Bob announces his opinion.

\begin{definition}\label{de1}
A \textit{test} is a function $T\dvtx \Delta\rightarrow\Sigma_{1}$.
\end{definition}

A test specifies an open set $T(P)$ considered inconsistent with an
opinion $P$. An expert who announces opinion $P$ is rejected on every
path $%
\omega$ belonging to $T ( P ) $. For the next definition, fix
$%
\varepsilon\in{}[0,1)$ and a subset $\Lambda$ of $\Delta$.

\begin{definition}\label{de2}
{A test} $\Lambda$\textit{-controls for type I error} {with
probability} $%
1-\varepsilon$ {if for any} $P\in\Lambda$,
\[
P \bigl( T ( P ) \bigr) \leq\varepsilon.
\]
\end{definition}

If a test $\Lambda$-controls for type I error, then an expert (with an
opinion in~$\Lambda$) expects to pass the test by honestly reporting what
he believes.

\subsection{Strategic forecasting}\label{sec2.1}

We now consider the case where Bob is uninformed about the odds of future
events, but may produce an opinion strategically in order to pass the test.
We allow strategic experts to select opinions at random. Let $\Delta
_{f}\Delta$ be the set of probability measures on $\Delta$ with finite
support. We call each $\zeta\in\Delta_{f}\Delta$ a
\emph{strategy}.\vadjust{\goodbreak}

\begin{definition}\label{de3}
{A test} {can be} \textit{manipulated} {with probability} $q\in [
0,1 ] $ {if there is a strategy }$\zeta$ such that for
every $\omega\in\Omega$,
\[
\zeta \bigl( \bigl\{ P\in\Delta\dvtx \omega\notin T ( P ) \bigr\} \bigr) \geq q.
\]
\end{definition}

If a test is manipulable with high probability, then a uninformed, but
strategic expert is likely to pass the test regardless of how the data
unfolds and how much data is available.

\begin{definition}\label{de4}
{A test} {is} \textit{nonmanipulable} {if for every strategy} $\zeta$
there is a cylinder $C_{\zeta}$ such that for every path $%
\omega\in C_{\zeta}$,
\[
\zeta \bigl( \bigl\{ P\in\Delta\dvtx \omega\notin T ( P ) \bigr\} \bigr) =0.
\]
\end{definition}

Nonmanipulable tests can reject uninformed experts. No matter which
strategy Bob employs, there is a finite history that, if observed,
discredits him. These are the only tests that are likely to pass informed
experts and may reject uninformed ones.

\section{Merging of opinions}\label{sec3}

We now review the main concepts behind the Blackwell--Dubins theorem.

\begin{definition}\label{de5}
{Let} $P,Q\in\mathbb{P}$. {The probability} $P$ \textit{merges} {with
}$Q$ if for every $\varepsilon>0$
\[
\lim_{t\rightarrow\infty}Q \Bigl( \Bigl\{ \omega\dvtx \sup
_{E\in\Sigma
} \bigl\llvert P \bigl( E|\omega^{t} \bigr) -Q
\bigl( E| \omega^{t} \bigr) \bigr\rrvert >\varepsilon \Bigr\} \Bigr) =0.
\]
\end{definition}

The expression $\sup_{E\in\Sigma}\llvert  P ( E|\omega
^{t} )
-Q ( E|\omega^{t} ) \rrvert $ is the distance between the
forecasts of $P$ and $Q$, conditional on the evidence available at time~$t$
and along the path~$\omega$. The probability $P$ merges with $Q$ if,
under $%
Q$, this distance goes to $0$ in probability. In particular, if $Q$
accurately describes the data generating process then the predictions
of $P$
are eventually accurate with high probability.

In this paper, merging is formulated in terms of convergence in probability
rather than almost sure convergence [as in \citet{BlaDub62}]. As
is well known, convergence in probability is particularly convenient in the
context of finitely additive probabilities. See, for instance, the
discussion in \citet{BerRig06}.

It is clear that for merging to occur, $P$ and $Q$ must be compatible
\emph{%
ex-ante}. The notion of absolute continuity formalizes this intuition.

\begin{definition}\label{de6}
Let $P,Q\in\mathbb{P}$. The probability $Q$ is
\textit{absolutely continuous} with respect to $P$, that is, $Q\ll P$,
if for every sequence of events $ ( E_{n} ) _{n=1}^{\infty
}$,
\[
\mbox{{if} }P ( E_{n} ) \rightarrow0\mbox{ {then} }%
Q (
E_{n} ) \rightarrow0.
\]
\end{definition}

If $P$ is $\sigma$-additive, then the definition is equivalent to requiring
that every event null under $P$ is also null under $Q$. Moreover,
if $P$ is a Savagean probability and $Q$ is a
probability satisfying $Q\ll P$, then $Q$ is Savagean as
well.

Absolute continuity is (essentially) necessary for merging.

\begin{proposition}\label{pr1}
Let $P,Q\in\mathbb{P}$ and $P ( \omega^{t} ) >0$ for every
cylinder $\omega^{t}$. If $P$ merges with $Q$ then $Q\ll P$.
\end{proposition}

In their seminal paper, Blackwell and Dubins show that when $P$ and $Q$
are $%
\sigma$-additive then absolute continuity suffices for merging.

\begin{theorem}[(Blackwell and Dubins)]\label{th1}
Let $P$ and $Q$ be $\sigma$-additive probability
measures on $ ( \Omega,\mathcal{B} ) $.  If $Q\ll P$, then
$P$ merges with $Q$.
\end{theorem}

One interpretation of the Blackwell--Dubins theorem is that multiple
repetitions of Bayes' rule lead to an agreement among agents who initially
hold compatible opinions. Another interpretation is that the
predictions of
Bayesian learners will eventually be accurate (provided that absolute
continuity holds). However, the Blackwell--Dubins theorem does not
extend to
all opinions. This motivates the next definition.

\begin{definition}\label{de7}
{A probability} $P\in\mathbb{P}$ {satisfies the} \textit{Blackwell--Dubins
property} if for every $Q\in\mathbb{P}$,
%
\begin{equation}
\mbox{{if} }Q\ll P\mbox{ {then} }P\mbox{ {merges with} }Q%
.
\end{equation}
Let $\Delta_{\mathrm{BD}}$ be the set of all opinions that satisfy the
Blackwell--Dubins property.
\end{definition}

So, an opinion $P$ satisfies the Blackwell--Dubins property if it
merges to
any compatible opinion $Q$. We show in the \hyperref[app]{Appendix} that $\Delta_{\mathrm{BD}}$ is
\emph{strictly} contained in the set of all opinions. That is, some opinions
satisfy the Blackwell--Dubins property, while others do not. We also
show that the set of probabilities satisfying the Blackwell--Dubins
property strictly contains the set of $\sigma$-additive probabilities.
We refer the reader to Example \ref{ex1}
and Theorem \ref{th10}, respectively.

Any exogenously given (or honestly reported) opinion can be tested and
rejected, whether or not the Blackwell--Dubins property holds. Thus, in the
absence of strategic considerations, the connection between the merging and
the testing of opinions is weak. We now show that this connection is much
stronger when there are incentive problems.

\section{Main results}\label{sec4}

\begin{theorem}\label{th2}
Consider the case where any opinion is allowed. Let $T$ be a test that $
\Delta$-controls for type I errors with probability $1-\varepsilon$. The test
$T$ can be manipulated with probability $1-\varepsilon-\delta$, for
every $%
\delta\in(0,1-\varepsilon]$.
\end{theorem}

\begin{theorem}\label{th3}
Consider the case where opinions must satisfy the Blackwell--Dubins property.
Fix $\varepsilon\in(0,1]$. There exists a test $T$ that $\Delta_{\mathrm{BD}}$-controls for type I error with probability $1-\varepsilon$ and is
nonmanipulable.
\end{theorem}

If Bob is free to announce any opinion, then he cannot be meaningfully tested
and discredited. Given any test that controls for type I error, Bob can
design a strategy which prevents rejection. However, if Bob is required to
announce opinions satisfying the Blackwell--Dubins property, then it is
possible to test and discredit him. These results show a strong connection
between the merging and the testing of opinions, but only when there are
incentive problems and agents may misrepresent what they know.

We now illustrate the basic ideas behind the proof of the two results. The
proof of Theorem \ref{th3} relies on a characterization of the set of probabilities
that satisfy the Blackwell--Dubins property. This characterization is also
crucial for the proof of Theorem \ref{th4} below. We show that $P\in\mathbb{P}$
satisfies the Blackwell--Dubins property if and only if it is an extreme
point of the set of probabilities $E ( P ) \subseteq\mathbb{P}$
which agree with $P$ on every cylinder.

The proof of necessity in this characterization is simple. Suppose, by
contradiction, that $P$ can be written as the convex combination
$P=\alpha
Q+ ( 1-\alpha ) R$, where $Q$ and $R$ belong to $E (
P ) $. Clearly, both $Q$ and $R$ are absolutely continuous with respect to $P$.
However, $P$ does not merge to $Q$ or $R$. The intuition is that $Q$
and $R$
agree on every finite history and so, the available data delivers equal
support to them.  The converse requires a deeper argument and relies
on Plachky's (\citeyear{Pla76}) theorem, which states that $P$ is an extreme
point of $E ( P ) $ if and only if the probability of every event
can be approximated by the probabilities of cylinders.

Given our characterization, the proof of Theorem \ref{th3} can be sketched as
follows: Let $P$ be an opinion satisfying the Blackwell--Dubins property.
Given that $P$ is strongly nonatomic, we can divide $\Omega$ into a
partition $ \{ A_{1},\ldots,A_{n} \} $ of events such that each of
them has probability less than $\varepsilon$. This property is a direct
implication
of Savage's postulate P6 and plays an important role in our result.
For general opinions, the events $ \{ A_{1},\ldots,A_{n} \} $
may have no useful structure and may not even be Borel sets. However,
since $P$ is an extreme point of $E ( P ) $, we can invoke
Plachky's theorem a second time and show that each $A_{i}$ can be
chosen to
be a cylinder. Now fix a path $\omega$. Let us say it belongs to
$A_{1}$. By
definition, there is a time $t$ such that $\omega^{t}=A_{1}$. We now define
a test $T$ such that $T ( P ) =\omega^{t}$ (note that the
period $%
t $ depends on the opinion $P$ because\vadjust{\goodbreak} the partition depends on it). By
definition, the test $\Delta_{\mathrm{BD}}$-controls for type I errors with
probability $1-\varepsilon$. Furthermore, it is a nonmanipulable test. Given
any strategy $\zeta$ we can find a period $m$ large enough such that $%
\omega^{m}$ rejects all opinions in the (finite) support of $\zeta$.
Therefore, in $\omega^{m}$, the probability of passing the test under $
\zeta$ is $0$.

We now sketch the proof of Theorem \ref{th2}. Consider a zero-sum game between
Nature and the expert. Nature chooses an opinion $P$ and the expert chooses
a strategy $\zeta$ (a random device producing opinions). The payoff of the
expert is the probability of passing the test. For each opinion $P$ chosen
by Nature there exists a strategy for the expert (to report $P$) that gives
him a payoff of at least $1-\varepsilon$. If Fan's (\citeyear{Fan53}) Minmax theorem
applies then there exists a strategy $\zeta  $that guarantees the expert
a payoff of at least $1-\varepsilon$ for \emph{every} opinion chosen by
Nature. In this case, the test is manipulable.

Fan's Minmax theorem requires Nature's action space to be compact and her
payoffs to be (lower semi) continuous. The main difficulty is that the set
of opinions is not compact in the natural topology, the weak* topology.
Hence, Fan's Minmax theorem cannot be directly applied. We consider a new
game, defined as above except that Nature can choose any probability in
$%
\mathbb{P}$ (not necessarily Savagean). By the Riesz representation and the
Banach--Alaoglu theorems, the set of all finitely additive probabilities
satisfy the necessary continuity and compactness conditions for Fan's Minmax
theorem. However, if now Nature chooses a \emph{non}-Savagean
probability $M$,
the expert cannot replicate her choice because he is restricted to opinions.

Based on the celebrated Hammer--Sobczyk decomposition theorem, we show the
following approximation result: For every $M\in \mathbb{P,}$ there is an
opinion $P$ such that $M ( U ) \leq P ( U ) $ for every
union $U$ of cylinders. Thus, $M ( T(P) ) \leq P (
T(P )
\leq\varepsilon$. It follows that if Natures chooses $M$ and the expert
chooses $P$ then he passes the test with probability at least
$1-\varepsilon$.
The proof is now concluded invoking Fan's Minmax theorem.

\section{Category tests}\label{sec5}

Theorem \ref{th3} provides conditions under which it is feasible to discredit
strategic experts. However, even a nonmanipulable test can be strategically
passed on some paths. Under $\sigma$-additivity, \citet{DekFei06}
and Olszewski and Sandroni (\citeyear{OlsSan09N1}) construct nonmanipulable \emph{category}
tests, where uninformed experts fail in all, but a topologically small (i.e.,
\emph{meager}) set of paths. We now show a difficulty in following this
approach in the general case of opinions that satisfy the Blackwell--Dubins
property.

\begin{definition}\label{de8}
A collection $\mathcal{I}$ of subsets of $\Omega$ is
a \textit{strictly proper ideal} if it satisfies the following properties:

\begin{longlist}[(1)]
\item[(1)] If $S\in\mathcal{I}$ and $R\subseteq S$ then $R\in\mathcal{I}$;

\item[(2)] If $R,S\in\mathcal{I}$ then $R\cup S\in\mathcal{I}$; and

\item[(3)] No cylinder belongs to $\mathcal{I}$.\vadjust{\goodbreak}
\end{longlist}
\end{definition}

A strictly proper ideal is a collection of sets which can be regarded as
``small.'' Property (1) is the natural
requirement that if a set $S$ is considered small then a set $R$ contained
in $S$ must also be considered small. Properties (2) and (3) are satisfied
by most commonly used notions of ``small''
sets, such as countable, nowhere dense, meager, sets of Lebesgue measure
zero and shy sets. To clarify our terminology, recall that an ideal is a
collection of subsets satisfying properties (1) and (2). An ideal is proper
if $\Omega$ does not belong to it. We refer to the elements of a strictly
proper ideal as \textit{small} sets and to their complements as \textit
{large} sets.

Strictly proper ideals can be defined in terms of probabilities.
Given $P\in\mathbb{P}$, define a set $N$ to be $P$-%
\emph{null} if there exists an event $E$ that satisfies $N\subseteq
E$ and $P ( E ) =0$. The collection of $P$-null sets is a
strictly proper ideal whenever $P$ satisfies $P ( \omega
^{t} ) >0$ for every cylinder $\omega^{t}$.

\begin{theorem}\label{th4}
Let $\mathcal{I}$ be a strictly proper ideal. There exists an opinion
$P\in
\Delta_{\mathrm{BD}}$ such that $P ( E ) =0$ for every event $E$ in $%
\mathcal{I}$.
\end{theorem}

There exists an opinion that satisfies the Blackwell--Dubins property and
finds all small events to be negligible. The proof of this result
relies on
the characterization of the set $\Delta_{\mathrm{BD}}$ of opinions satisfying the
Blackwell--Dubins property that we discussed in the previous section. Theorem
\ref{th4} shows a basic tension between the control of type I errors and
the use of genericity arguments. Suppose Alice intends to design a test that
discredits Bob on a large set of paths, irrespectively of his strategy.
Then the set of paths $ ( T(P) ) ^{c}$ that do not reject
opinion $%
P$ must be small [otherwise Bob could simply announce opinion $P$ and pass
the test on $(T(P))^{c}$, a~nonsmall set of realizations]. But if $P$ is
the opinion obtained from Theorem \ref{th4}, we must have $P (
(T(P))^{c} )
=0$. So, the test cannot control for type I errors. We have just proved the
following corollary.

\begin{corollary}\label{co1}
Let $\mathcal{I}$ be a strictly proper ideal. For every test $T$ which $
\Delta_{\mathrm{BD}}$-controls type I errors with positive probability there exists
a strategy $\zeta$ such that the set
\[
\bigl\{ \omega\dvtx \zeta \bigl( \bigl\{ P\dvtx \omega\notin T ( P ) \bigr\}
\bigr) =1 \bigr\}
\]
is not small.
\end{corollary}

Thus, the stronger nonmanipulable tests in \citet{DekFei06} and
Olszewski and Sandroni (\citeyear{OlsSan09N1}) cannot be obtained in the general case of
opinions that satisfy the Blackwell--Dubins property.

\section{Extensions}\label{sec6}

In this section, we extend our analysis to the case where the set of
per-period outcomes may be infinite.

\subsection{Setup}

Let $\mathcal{X}$ be a separable metric space of outcomes and denote by
$%
\Omega$ the set of paths $\mathcal{X}^{\infty}$. As before, $\omega^{t}$
is the cylinder of length $t\geq0$ with base $\omega\in\Omega$ (in
particular $\omega^{0}=\Omega$). The set $\Omega$ is endowed with the
product topology and a $\sigma$-algebra $\Sigma$ containing all open sets.
We denote by $\mathbb{P}$ the set of finitely additive probabilities on
$%
( \Omega,\Sigma ) $ and by $\Delta$ the subset of opinions
(i.e., strongly nonatomic probabilities).

\subsection{Conditional probabilities}\label{sec6.2}

Let $\mathcal{H}$ be the set of all cylinders. We say that a function
\[
P\dvtx \Sigma\times\mathcal{H\rightarrow} [ 0,1 ]
\]
is a \emph{conditional probability} if for every $t\geq0$ and $\omega
\in
\Omega$:

\begin{longlist}[(1)]
\item[(1)]$P ( \cdot |\omega^{t} ) $ $\in\mathbb{P}$;

\item[(2)]$P ( \omega^{t}|\omega^{t} ) =1$; and

\item[(3)]$P ( E\cap\omega^{t+n}|\omega^{t} ) =P ( E|\omega
^{t+n} ) P ( \omega^{t+n}|\omega^{t} ) $ for any event $E$
and $n\geq0$.
\end{longlist}

The definition of conditional probability follows \citet{BerRegRig97}, where properties (1)--(3) are justified on the basis of de Finetti's
\emph{coherence} principle: a real function $P$ defined on $\Sigma
\times
\mathcal{H}$ satisfies properties (1)--(3) if and only if a bookie, who
sets $%
P ( E | \omega^{t} ) $ as the price of a conditional bet on
event~$E$, cannot incur in a Dutch book. We refer the reader to Regazzini
(\citeyear{Reg85,Reg87}), \citet{deF90} and \citet{BerRig02} for a precise
statement and a formal discussion.

A conditional probability is a \emph{conditional opinion} if $P (
\cdot | \Omega ) $ is strongly nonatomic. We denote by $\mathbb
{P%
}^{\ast}$ and by $\Delta^{\ast}$ the sets of conditional probabilities
and conditional opinions, respectively. To simplify the exposition,
given an
event $E$ and a conditional probability $P$, we use the notation $P (
E ) $ instead of the more precise $P ( E | \Omega ) $.

At time $0$, Bob is required to announce a conditional opinion $P$. So, a
test $T$ is now a function mapping each conditional opinion $P$ to an open
subset $T ( P ) $ of $\Omega$. The definitions of type~I errors,
manipulable and nonmanipulable tests are analogous to the definitions of
Section~\ref{sec2} and can be obtained by replacing $\Delta$ with~$\Delta
^{\ast}$.

\subsection{Merging}\label{sec6.3}

We now extend the definition of merging of opinions. We say that
$\mathcal{X}
$ is a \emph{discrete space} if it is countable and endowed with the
discrete topology.

\begin{definition}\label{de9}
Let $\mathcal{X}$ be a discrete space. If $P,Q\in
\mathbb{P}^{\ast}$, the conditional probability $P$ \textit{merges}
with $Q$ if for every $\varepsilon>0$,
\[
\lim_{t\rightarrow\infty}Q \Bigl( \Bigl\{ \omega\dvtx \sup
_{E\in\Sigma
} \bigl\llvert P \bigl( E|\omega^{t} \bigr) -Q
\bigl( E| \omega^{t} \bigr) \bigr\rrvert >\varepsilon \Bigr\} \Bigr) =0.
\]
\end{definition}

The next definition is based on \citet{BlaDub62}.\vadjust{\goodbreak}

\begin{definition}\label{de10}
Let $\mathcal{X}$ be a discrete space. A conditional
probability $P$ satisfies the \textit{Blackwell--Dubins property} if for
every probability $Q\in\mathbb{P}$ such that $Q\ll P ( \cdot
| \Omega ) $ there exists a conditional probability $%
\widetilde{Q}$ such that
\[
\widetilde{Q} ( \cdot | \Omega ) =Q\mbox{ {and} }P\mbox{ {merges with} }
\widetilde{Q}.
\]
Let $\Delta_{\mathrm{BD}}^{\ast}$ be the set of all conditional
opinions that satisfy the Blackwell--Dubins property.
\end{definition}

We now show that the connection between testability and merging of opinions
extends to this setup.

\subsection{Results}\label{sec6.4}

\begin{theorem}\label{th5}
Let $\mathcal{X}$ be a separable metric space.
Consider the case where any conditional opinion is allowed. Let $T$ be a
test that $\Delta^{\ast}$-controls for type I errors with probability
$%
1-\varepsilon$. The test $T$ can be manipulated with probability
$1-\varepsilon
-\delta$, for every $\delta\in(0,1-\varepsilon]$.
\end{theorem}

\begin{theorem}\label{th6}
Let $\mathcal{X}$ be a discrete space. Consider the case where conditional
opinions must satisfy the Blackwell--Dubins property. Fix $\varepsilon\in
(0,1]$%
. There exists a test $T$ that $\Delta_{\mathrm{BD}}^{\ast}$-controls for type I
error with probability $1-\varepsilon$ and is nonmanipulable.
\end{theorem}

If it is possible for Bob to announce any conditional opinion, then he cannot
be meaningfully tested and discredited. If Bob is restricted to conditional
opinions satisfying the Blackwell--Dubins property, then it is possible to
test and discredit him.

The proof of Theorem \ref{th5} follows the proof of Theorem \ref{th2}. The proof of Theorem
\ref{th6} is based on the following result: a conditional opinion $P\in\Delta
_{\mathrm{BD}}^{\ast}$ satisfies $\lim_{t\rightarrow\infty}P ( \omega
^{t} ) =0$ for every path $\omega$. This step requires a new argument,
because the characterization of Blackwell--Dubins property used in the proof
of Theorem \ref{th3} does not readily extend to the case where $\mathcal{X}$ is
infinite. Once this continuity property is shown to hold, the proof
continues as in Theorem \ref{th3}. We fix a path $\omega$ and for each $P\in
\Delta
_{\mathrm{BD}}^{\ast}$ we choose a large enough period $t_{P}$ such that $P (
\omega^{t_{P}} ) <\varepsilon$. Because $\mathcal{X}$ is assumed to
be a
discrete space, each cylinder $\omega^{t_{P}}$ is open. Therefore, we can
define test $T$ such that $T ( P ) =\omega^{t_{P}}$ for every
$%
P\in\Delta_{\mathrm{BD}}^{\ast}$. Following the proof of Theorem \ref{th3}, we show
that $T$
is nonmanipulable.

\begin{appendix}\label{app}
\section{}

We now provide an example of an opinion that violates the Blackwell--Dubins
property.

\begin{example}\label{ex1}
Let $\Omega= \{ 0,1 \} ^{\infty}$ and $\Sigma=\mathcal{B}$.
Denote by $X_{1},X_{2},\ldots $ the coordinate projections on $\Omega$. For
every $n\geq1$, let $P_{n}$ be the $\sigma$-additive probability
defined as
\[
P_{n} ( X_{k}=0 ) =2^{-k}\qquad\mbox{for }k\leq
n\quad \mbox{and}\quad P_{n} ( X_{k}=0 ) =1\qquad\mbox{for
}k>n
\]
and let $P_{\infty}$ be the $\sigma$-additive probability defined as $
P_{\infty} ( X_{k}=0 ) =2^{-k}$ for all $k$.

Consider the opinion $P=\frac{1}{2}P_{\infty}+\frac{1}{2}\int
P_{n}\,d\lambda
( n ) $, where $\lambda$ is a finitely additive probability
on $%
(
\mathbb{N},2^{%
\mathbb{N}
} ) $ such that $\lambda (  \{ n \}  ) =0$ for every
$n$. The finite additivity of the mixture $\int P_{n}\,d\lambda (
n ) $ may reflect the difficulty of predicting \emph{when} the
per-period probability of observing the outcome $0$ will change from $0.5$
to $1$.

Clearly, $P_{\infty}\ll P$. However, $P$ does not merge with
$P_{\infty}$.
To this end, let $A$ be the set of all paths where the outcome $1$ appears
infinitely often. Then $P_{n} ( A ) =0$ for every $n$ and $%
P_{\infty} ( A ) =1$. For every cylinder $\omega^{t}$, we
have
\[
P \bigl( \omega^{t} \bigr) =\frac{1}{2}P_{\infty} \bigl(
\omega ^{t} \bigr) +%
\frac{1}{2}\int
_{ \{ n:n>t \} }P_{n} \bigl( \omega^{t} \bigr) \,d
\lambda ( n ) =P_{\infty} \bigl( \omega^{t} \bigr)
\]
and moreover,
\begin{eqnarray*}
P_{\infty} \bigl( A|\omega^{t} \bigr) -P \bigl( A|
\omega^{t} \bigr) &=&1-%
\frac{({1}/{2})P_{\infty} ( A\cap\omega^{t} ) +({1}/{2})
\int_{n}P_{n} ( A\cap\omega^{t} ) \,d\lambda ( n ) }{
P_{\infty} ( \omega^{t} ) }
\\
&=&1-\frac{1}{2}P_{\infty} \bigl( A|\omega^{t} \bigr)
\\
&=&\frac{1}{2}
\end{eqnarray*}
for every $\omega$ and every $t$. Thus, $P$ does not merge with
$P_{\infty}
$.
\end{example}

\subsection{Preliminaries}

To minimize repetitions, throughout the \hyperref[app]{Appendix} $\Omega$ stands
for either $ \{ 0,1 \} ^{\infty}$ or $\mathcal{X}^{\infty}$. For
every algebra $\mathcal{A}$ of subsets of $\Omega$ denote by $\mathbb
{P}%
( \mathcal{A} ) $ the space of finitely additive probabilities
defined on $ ( \Omega,\mathcal{A} ) $. When $\mathcal
{A}=\Sigma$,
we write $\mathbb{P}$ instead of $\mathbb{P} ( \Sigma ) $. We
denote by $\Delta\subseteq\mathbb{P}$ the set of opinions (strongly
nonatomic probabilities). The space $\mathbb{P} ( \mathcal{A}
) $
is endowed with the weak* topology. It is the coarsest topology for which
the functional $P\mapsto\int\varphi \,dP$ is continuous for every
function $%
\varphi\dvtx \Omega\rightarrow
\mathbb{R}
$ that has finite range and is measurable with respect to~$\mathcal{A}$.
This should not be confused with the more common weak* topology
generated by
bounded \emph{continuous} functions.

\section{Merging of opinions}\label{sec8}

In this subsection, we describe of the set of opinions on $\Omega=
\{
0,1 \} ^{\infty}$ that satisfy the Blackwell--Dubins property. We first
show that absolute continuity is essentially a necessary condition for
merging.

\begin{pf*}{Proof of Proposition \ref{pr1}}
Assume $Q$ is not absolutely continuous with respect to $P$. Then there
exists a sequence of events $ ( E_{n} ) _{n=1}^{\infty}$ and
some $%
\alpha>0$ such that $P ( E_{n} ) \rightarrow0$ but $Q (
E_{n} ) >\alpha$ for every $n$. Suppose $P$ merges with $Q$. For every
$t$, let $\mathcal{C}_{t}$ be a collection of pairwise disjoint
cylinders of
length $t$ such that $Q ( \omega^{t} ) >0$ for every $\omega
^{t}\in\mathcal{C}_{t}$ and $Q ( \cup\mathcal{C}_{t} ) =1$.
Fix $\delta\in ( 0,\frac{\alpha}{4} ) $. There exists a
time $T$
large enough such that for every $t\geq T$ there is a subset $\mathcal
{D}%
_{t}\subseteq\mathcal{C}_{t}$ such that $Q ( \cup\mathcal{D}%
_{t} ) \geq1-\delta$  and $\sup_{E\in\Sigma}\llvert  Q (
E|\omega^{t} ) -P ( E|\omega^{t} ) \rrvert \leq
\delta$
for every $\omega^{t}\in\mathcal{D}_{t}$. Because $P ( \omega
^{t} ) >0$ for every $\omega^{t}$, the expression $P (
E|\omega
^{t} ) $ is well defined. For every event $E_{n}$,
\[
Q ( E_{n} ) =\sum_{\omega^{t}\in\mathcal{D}_{t}}Q \bigl(
E_{n}|\omega^{t} \bigr) Q \bigl( \omega^{t} \bigr)
+\sum_{\omega
^{t}\in
\mathcal{C}_{t}-\mathcal{D}_{t}}P \bigl( E|\omega^{t} \bigr) P
\bigl( \omega ^{t} \bigr)
\]
hence, $\sum_{\omega^{t}\in\mathcal{D}_{t}}Q ( E_{n}|\omega
^{t} )
Q ( \omega^{t} ) \geq\alpha-\delta$. Define
\[
\mathcal{E}_{n}= \biggl\{ \omega^{t}\in
\mathcal{D}_{t}\dvtx Q \bigl( E_{n}|\omega^{t}
\bigr) \geq\frac{\alpha}{2} \biggr\}.
\]
We have
\begin{eqnarray*}
\alpha-\delta&\leq&\sum_{\omega^{t}\in\mathcal{D}_{t}}Q \bigl(
E_{n}|\omega^{t} \bigr) Q \bigl( \omega^{t} \bigr)
\\
&=&\sum_{\omega^{t}\in\mathcal{E}_{n}}Q \bigl( E_{n}|
\omega^{t} \bigr) Q \bigl( \omega^{t} \bigr) +\sum
_{\omega^{t}\in\mathcal
{D}_{t}-\mathcal{E}%
_{n}}Q \bigl( E_{n}|\omega^{t} \bigr) Q
\bigl( \omega^{t} \bigr)
\\
&\leq&Q ( \cup\mathcal{E}_{n} ) + \frac{\alpha}{2}Q ( \cup
\mathcal{D}_{t}-\cup \mathcal{E}_{n} )
\\
&\leq&Q ( \cup\mathcal{E}_{n} ) + \frac{\alpha}{2}
\end{eqnarray*}
hence, $Q ( \cup\mathcal{E}_{n} ) \geq\frac{\alpha
}{2}-\delta$%
. Now let $n^{\ast}$ be large enough such that $\frac{\alpha
}{2}-P (
E_{n^{\ast}}|\omega^{t} ) >\delta$ for every $\omega$. Then, for
every $\omega^{t}\in\mathcal{E}_{n^{\ast}}$,
\begin{eqnarray*}
\sup_{E\in\Sigma} \bigl\llvert Q \bigl( E|\omega^{t} \bigr)
-P \bigl( E|\omega ^{t} \bigr) \bigr\rrvert &\geq&Q \bigl(
E_{n^{\ast}}|\omega^{t} \bigr) -P \bigl( E_{n^{\ast}}|
\omega^{t} \bigr)
\\
&\geq&\frac{\alpha}{2}-P \bigl( E_{n^{\ast}}|\omega^{t} \bigr)
\\
&>&\delta.
\end{eqnarray*}
To summarize, $Q (  \{ \omega\dvtx \sup_{E\in\Sigma}\llvert
Q (
E|\omega^{t} ) -P ( E|\omega^{t} ) \rrvert >\delta
\}  ) >\frac{\alpha}{2}-\delta$ for every $t\geq T$.
Therefore, $%
Q ( \cup\mathcal{D}_{t} ) \leq1- ( \frac{\alpha}{2}%
-\delta ) $. By definition, $Q ( \cup\mathcal{D}_{t} )
\geq1-\delta$; hence, $\delta\geq\frac{\alpha}{4}$. A contradiction.
Hence, $P$ does not merge with $Q$.
\end{pf*}

Our main result is a characterization of the set of opinions that satisfy
the Blackwell--Dubins property. We first recall some results on
extensions of
finitely additive probabilities. Let $\mathcal{A}_{1}$ and $\mathcal{A}_{2}$
be two algebras of subsets of $\Omega$ such that $\mathcal
{A}_{1}\subseteq
\mathcal{A}_{2}$. Given $P\in\mathbb{P} ( \mathcal{A}_{1} ) $
and $%
Q\in\mathbb{P} ( \mathcal{A}_{2} ) $, call $Q$ an \emph{extension
of} $P$ \emph{from} $\mathcal{A}_{1}$ \emph{to} $\mathcal{A}_{2}$ if $%
P ( A ) =Q ( A ) $ for every $A\in\mathcal{A}_{1}$.
Let $%
E ( P,\mathcal{A}_{1},\mathcal{A}_{2} ) $ be the set of extensions
of $P$ from $\mathcal{A}_{1}$ to $\mathcal{A}_{2}$. As is well known, the
set $E ( P,\mathcal{A}_{1},\mathcal{A}_{2} ) $ is nonempty.
Moreover, it is a convex and compact subset of $\mathbb{P} (
\mathcal{A}%
_{2} ) $. The set of extreme points of $E ( P,\mathcal{A}_{1},%
\mathcal{A}_{2} ) $ has been studied in great generality. We refer the
reader to \citet{Lip07} and \citet{Pla76} for further results and
references.

\begin{theorem}[{[\citet{Pla76}]}]\label{th7}
Fix two algebras $\mathcal{A}_{1} \subseteq\mathcal{A}_{2}$ and $P\in
\mathbb{P} ( \mathcal{A}_{1} ) $. A~probability $Q\in E (
P,%
\mathcal{A}_{1},\mathcal{A}_{2} ) $ is an extreme point of $E
( P,%
\mathcal{A}_{1},\mathcal{A}_{2} ) $ if and only if for every
$\varepsilon
>0$ and $A_{2}\in\mathcal{A}_{2}$ there exists $A_{1}\in\mathcal{A}_{1}$
such that $Q ( A_{2}\bigtriangleup A_{1} ) <\varepsilon$.
\end{theorem}

Let $\mathcal{F}$ be the algebra generated by all cylinders of $ \{
0,1 \} ^{\infty}$. An event belongs to $\mathcal{F}$ if and only
if it
is a finite union of (pairwise disjoint) cylinders. Recall that
$\mathcal{B}$
is the Borel $\sigma$-algebra induced on $\Omega$ by the product topology.
Then $\mathcal{F}\subseteq\mathcal{B}\subseteq\Sigma$. For every
$P\in
\mathbb{P}$ let $P_{\mathcal{F}}$ be the restriction of $P$ on $\mathcal
{F}$%
. It is easy to see that $P_{\mathcal{F}}$ is $\sigma$-additive. By
Carath%
\'{e}odory theorem, it admits a $\sigma$-additive extension from
$\mathcal{F}
$ to~$\mathcal{B}$, denoted by $P_{\sigma}$.

The next result is well known.

\begin{lemma}\label{le1}
Let $\Omega= \{ 0,1 \} ^{\infty}$. For all $Q$ and $P$ in $%
\mathbb{P}$, if $Q_{\mathcal{F}}\ll P_{\mathcal{F}}$ then $Q_{\sigma
}\ll
P_{\sigma}$.
\end{lemma}

We can now state our main result on merging.

\begin{theorem}\label{th8}
Let $\Omega= \{ 0,1 \} ^{\infty}$. For every $P\in\mathbb{P,}$
the following are equivalent:

\begin{longlist}[(1)]
\item[(1)] $P$ is an extreme point of $E ( P_{\mathcal{F}},\mathcal
{F},\Sigma
) $.

\item[(2)] $P$ satisfies the Blackwell--Dubins property.

\item[(3)] For all $Q,R\in\mathbb{P}$, if $P=\alpha Q+ ( 1-\alpha
) R$
for some $\alpha\in ( 0,1 ) $ then $P$ merges with $Q$ and $R$.
\end{longlist}
\end{theorem}

\begin{pf}
(1)${}\implies{}$(2). Let $P$ be an extreme point of $E (
P_{\mathcal{F}%
},\mathcal{F},\Sigma ) $. If $Q\ll P$, then $Q_{\sigma}\ll
P_{\sigma}$
by Lemma \ref{le1}. By the Blackwell--Dubins theorem,
\[
Q_{\sigma} \Bigl( \Bigl\{ \omega\dvtx \lim_{t\rightarrow\infty} \Bigl(
\sup_{B\in\mathcal{B}} \bigl\llvert Q_{\sigma} \bigl( B|
\omega^{t} \bigr) -P_{\sigma} \bigl( B|\omega^{t} \bigr)
\bigr\rrvert \Bigr) =0 \Bigr\} \Bigr) =1.
\]
In particular, the sequence of random variables
\[
\Bigl( \omega\mapsto\sup_{F\in\mathcal{F}} \bigl\llvert Q_{\sigma}
\bigl( F|\omega^{t} \bigr) -P_{\sigma} \bigl( F|
\omega^{t} \bigr) \bigr\rrvert \Bigr) _{t=1}^{\infty}
\]
converges to $0$, $Q_{\sigma}$-almost surely. Therefore, the sequence
converges in probability. For each $\varepsilon>0$,
\[
\lim_{t\rightarrow\infty}Q_{\sigma} \Bigl( \Bigl\{ \omega\dvtx \sup
_{F\in
\mathcal{F}} \bigl\llvert Q_{\sigma} \bigl( F|
\omega^{t} \bigr) -P_{\sigma
} \bigl( F|\omega^{t} \bigr)
\bigr\rrvert >\varepsilon \Bigr\} \Bigr) =0.
\]
Since the last expression only involves events belonging to $\mathcal
{F}$,
\[
\lim_{t\rightarrow\infty}Q \Bigl( \Bigl\{ \omega\dvtx \sup
_{F\in\mathcal
{F}%
} \bigl\llvert Q \bigl( F|\omega^{t} \bigr) -P
\bigl( F| \omega^{t} \bigr) \bigr\rrvert >\varepsilon \Bigr\} \Bigr) =0.
\]

The proof is complete by showing that
\[
\sup_{E\in\Sigma} \bigl\llvert Q \bigl( E|\omega^{t} \bigr)
-P \bigl( E|\omega ^{t} \bigr) \bigr\rrvert =\sup_{F\in\mathcal{F}}
\bigl\llvert Q \bigl( F|\omega ^{t} \bigr) -P \bigl( F|
\omega^{t} \bigr) \bigr\rrvert
\]
for every $\omega^{t}$ such that $Q ( \omega^{t} ) >0$.

To this end, fix an event $E\in\Sigma$ and a cylinder $\omega^{t}$ such
that $Q ( \omega^{t} ) >0$. By Plachky's theorem, there
exists a
sequence of events $ ( F_{n} ) _{n=1}^{\infty}$ in $\mathcal{F}$
such that $P ( E\bigtriangleup F_{n} ) \rightarrow0$ as $%
n\rightarrow\infty$. For each $n$, the inequality
\begin{eqnarray*}
\bigl\llvert Q \bigl( E|\omega^{t} \bigr) -P \bigl( E|
\omega^{t} \bigr) \bigr\rrvert &\leq& \bigl\llvert Q \bigl( E|
\omega^{t} \bigr) -Q \bigl( F_{n}|\omega^{t} \bigr)
\bigr\rrvert
\\
&&{}+ \bigl\llvert Q \bigl( F_{n}|\omega ^{t} \bigr) -P
\bigl( F_{n}|\omega^{t} \bigr) \bigr\rrvert
\\
&&{}+ \bigl\llvert P \bigl( F_{n}|\omega^{t} \bigr) -P \bigl(
E|\omega ^{t} \bigr) \bigr\rrvert
\end{eqnarray*}
implies
\begin{eqnarray*}
&&\bigl\llvert Q \bigl( E|\omega^{t} \bigr) -P \bigl( E|
\omega^{t} \bigr) \bigr\rrvert \\
&&\qquad\leq Q \bigl( E\bigtriangleup
F_{n}|\omega^{t} \bigr) +\sup_{F\in\mathcal{F}} \bigl
\llvert Q \bigl( F|\omega^{t} \bigr) -P \bigl( F|\omega^{t}
\bigr) \bigr\rrvert +P \bigl( E\bigtriangleup F_{n}|\omega
^{t} \bigr).
\end{eqnarray*}
Because $P ( E\bigtriangleup F_{n} ) \rightarrow0$ and $Q (
\omega^{t} ) >0$, it follows that $P ( \omega^{t} ) >0$
and $P (
E\bigtriangleup F_{n}|\omega^{t} ) \rightarrow0$. Absolute continuity
implies $Q ( E\bigtriangleup F_{n}|\omega^{t} ) \rightarrow0$.
Therefore,
\[
\bigl\llvert Q \bigl( E|\omega^{t} \bigr) -P \bigl( E|
\omega^{t} \bigr) \bigr\rrvert \leq\sup_{F\in\mathcal{F}} \bigl
\llvert Q \bigl( F|\omega ^{t} \bigr) -P \bigl( F|\omega^{t}
\bigr) \bigr\rrvert
\]
thus
\[
\sup_{E\in\Sigma} \bigl\llvert Q \bigl( E|\omega^{t} \bigr)
-P \bigl( E|\omega ^{t} \bigr) \bigr\rrvert \leq\sup
_{F\in\mathcal{F}} \bigl\llvert Q \bigl( F|\omega^{t} \bigr) -P
\bigl( F|\omega^{t} \bigr) \bigr\rrvert
\]
as claimed.

(2)${}\implies{}$(3). If $P=\alpha Q+ ( 1-\alpha ) R$ then $Q\ll P$,
hence $P$ merges with $Q$.

(3)${}\implies{}$(1). Assume by way of contradiction that $P$ is not an extreme
point of $E ( P_{\mathcal{F}},\mathcal{F},\Sigma ) $. Then there
exist $Q,R$ in $E ( P_{\mathcal{F}},\mathcal{F},\Sigma ) $ such
that$  P=\alpha Q+ ( 1-\alpha ) R $, $\alpha\in (
0,1 ) $ and $Q\neq R$. By assumption, $P$ merges with $Q$. Let $%
\mathcal{C}_{t}$ be a collection of pairwise disjoint cylinders of
length $t$
such that $Q ( \omega^{t} ) >0$ for every $\omega^{t}\in
\mathcal{%
C}_{t}$ and $Q ( \cup\mathcal{C}_{t} ) =1$. For every $t$ large
enough, there exists a subset $\mathcal{D}_{t}\subseteq\mathcal
{C}_{t}$ such
that $Q ( \cup\mathcal{D}_{t} ) \geq1-\varepsilon$  and $%
\sup_{E\in\Sigma}\llvert  Q ( E|\omega^{t} ) -P (
E|\omega
^{t} ) \rrvert \leq\varepsilon$ for every $\omega^{t}\in
\mathcal{D}%
_{t}$. For every event $E$,
\begin{eqnarray*}
\bigl\llvert Q ( E ) -P ( E ) \bigr\rrvert &=& \biggl\llvert \sum
_{\omega^{t}\in\mathcal{C}_{t}}Q \bigl( E|\omega^{t} \bigr) Q \bigl(
\omega^{t} \bigr) -\sum_{\omega^{t}\in\mathcal{C}_{t}}P \bigl( E|
\omega ^{t} \bigr) P \bigl( \omega^{t} \bigr) \biggr\rrvert
\\
&=& \biggl\llvert \sum_{\omega^{t}\in\mathcal{C}_{t}} \bigl( P \bigl( E|
\omega ^{t} \bigr) -Q \bigl( E|\omega^{t} \bigr) \bigr) Q
\bigl( \omega ^{t} \bigr) \biggr\rrvert
\\
&\leq& \sum_{\omega\in\mathcal{D}_{t}} \bigl\llvert P \bigl( E|\omega
^{t} \bigr) -Q \bigl( E|\omega^{t} \bigr) \bigr\rrvert Q
\bigl( \omega ^{t} \bigr)
\\
&&{}+ \sum_{\omega\in\mathcal{C}_{t}-\mathcal{D}_{t}} \bigl\llvert P \bigl( E|
\omega^{t} \bigr) -Q \bigl( E|\omega^{t} \bigr) \bigr\rrvert Q
\bigl( \omega^{t} \bigr)
\\
&\leq& \varepsilon Q ( \cup\mathcal{D}_{t} ) + \bigl( 1-Q ( \cup
\mathcal{D}_{t} ) \bigr)
\\
&\leq&2\varepsilon,
\end{eqnarray*}
where the first two equalities follow from the fact that $Q ( \omega
^{t} ) =P ( \omega^{t} ) $ for all $\omega^{t}$. Since $E$
and $\varepsilon$ are arbitrary, we have $P=Q$. A contradiction.
Therefore, $P$
must be an extreme point of $E ( P_{\mathcal{F}},\mathcal{F},\Sigma
) $.
\end{pf}

The next result shows a useful property of opinions that satisfy the
Blackwell--Dubins property.

\begin{theorem}\label{th9}
Let $\Omega= \{ 0,1 \} ^{\infty}$. For every $P\in\Delta_{\mathrm{BD}}$
and every $\varepsilon>0$, there exists a partition $ \{
C_{1},\ldots,C_{n} \} $ of $\Omega$ such that for each $i=1,\ldots,n$, $
C_{i} $ is a cylinder and $P ( C_{i} ) \leq\varepsilon$.
\end{theorem}

\begin{pf}
Let $P\in\Delta_{\mathrm{BD}}$ and fix $\varepsilon>0$. Because $P$ is strongly
nonatomic, then there exists a partition $ \{ E_{1},\ldots,E_{m}
\} $
of events such that $P ( E_{i} ) <\frac{\varepsilon}{2}$ for
every $%
i=1,\ldots,m$. By Theorem \ref{th8}, $P$ is an extreme point of $E (
P_{\mathcal{F}%
},\mathcal{F},\Sigma ) $. By Plachky's theorem, for each $i$ we can
find a sequence $ ( F_{i,k} ) _{k=1}^{\infty}$ in $\mathcal{F}$
such that $P ( E_{i}\bigtriangleup F_{i,k} ) \rightarrow0$ as
$%
k\rightarrow\infty$. Choose $K$ large enough such that $P (
E_{i}\bigtriangleup F_{i,K} ) <\frac{\varepsilon}{2m}$ for each $i$.

Let $F_{1}=F_{1,K}$ and define $F_{i}=F_{i,K}-\bigcup_{j=1}^{i-1}F_{j,K}$ for
each $i=2,\ldots,m$. Let $F_{m+1}=\Omega-\bigcup_{i=1}^{m}F_{i,K}$ and consider
the partition $ \{ F_{1},\ldots,F_{m+1} \} $. It satisfies $P (
F_{i} ) \leq P ( F_{i,K} ) <\frac{\varepsilon}{2}+\frac
{\varepsilon
}{2m}$ for each $i=2,\ldots,m$. Moreover,
\[
P ( F_{m+1} ) =P \Biggl( \Biggl( \bigcup_{i=1}^{m}E_{i}
\Biggr) - \Biggl( \bigcup_{i=1}^{m}F_{i,K}
\Biggr) \Biggr) \leq P \Biggl( \bigcup_{i=1}^{m}
( E_{i}-F_{i,K} ) \Biggr) \leq\frac{\varepsilon}{2}.
\]
Therefore, $P ( F_{i} ) \leq\varepsilon$ for each $F_{i}\in
\{
F_{1},\ldots,F_{m+1} \} $. Because each $F_{i}$ is a finite union of
pairwise disjoint cylinders, the proof is complete.
\end{pf}

The next theorem shows that for every sigma additive probability $P$ we can
find a continuum of probabilities that agree with $P$ on every cylinder,
fail $\sigma$-additivity but satisfy the Blackwell--Dubins property. A
related result appears in \citet{Lip01}.

\begin{theorem}\label{th10}
Let $\Omega= \{ 0,1 \} ^{\infty}$. For every $\sigma$-additive
probability $P\in\mathbb{P,}$ the set
\[
\{ Q\in\Delta_{\mathrm{BD}}\dvtx Q_{\mathcal{F}}=P_{\mathcal{F}}, Q\mbox{ is
not } \sigma\mbox{-additive} \}
\]
has cardinality at least $\mathfrak{c}$.
\end{theorem}

\begin{pf}
Fix a collection $ \{D_{\xi}\dvtx \xi\in [ 0,1 ]  \} $ of
pairwise disjoint, Borel and dense subsets of $\Omega$ [see, e.g.,
\citet{Ced66}]. For every $\xi\in [ 0,1 ] $, let $\mathcal
{A}_{\xi}
$ be the algebra generated by $\mathcal{F}\cup \{ D_{\xi} \} $.
That is,
\[
A_{\xi}= \bigl\{ ( F_{1}\cap D_{\xi} ) \cup \bigl(
F_{2}\cap D_{\xi}^{c} \bigr) \dvtx
F_{1},F_{2} \in\mathcal{F} \bigr\}.
\]
Let $\rho_{\xi}\in\mathbb{P} ( \mathcal{A}_{\xi} ) $ be defined
as
\[
\rho_{\xi} ( F\cap D_{\xi} ) =P ( F )
\]
for every $F\in\mathcal{F}$. Because $D_{\xi}$ is dense, then $F\cap
D_{\xi
}\neq\varnothing$ for every $F\in\mathcal{F}$. Hence, $\rho_{\xi}$
is well
defined. It satisfies $\rho_{\xi}  ( D_{\xi} ) =1$.

For each $\xi\in [ 0,1] $ fix an extreme point $P_{\xi}$
of $%
E ( \rho_{\xi},\mathcal{A}_{\xi},\Sigma ) $. We claim it is
also an extreme point of $E ( P_{\mathcal{F}},\mathcal{F},\Sigma
)
$. By construction, $P_{\xi}\in E ( P_{\mathcal{F}},\mathcal
{F},\Sigma
) $. Now suppose $P_{\xi}=\alpha Q+(1-\alpha)R$, with $\alpha
\in%
[ 0,1 ] $ and $Q,R\in E ( P_{\mathcal{F}},\mathcal
{F},\Sigma
) $. Because $P_{\xi} ( D_{\xi} ) =1$, then $Q (
D_{\xi
} ) =R ( D_{\xi} ) =1$. Hence, $P ( F ) =Q (
F\cap D_{\xi} ) =R ( F\cap D_{\xi} ) $ for every $F\in
\mathcal{F}$. Therefore, $Q,R\in E ( \rho_{\xi},\mathcal{A}_{\xi
},\Sigma ) $. By assumption, $P_{\xi}$ is an extreme point of $%
E ( \rho_{\xi},\mathcal{A}_{\xi},\Sigma ) $. Hence, $P=Q=R$.
This concludes the proof of the claim.

By Theorem \ref{th8}, each $P_{\xi}$ satisfies the Blackwell--Dubins property.
Moreover, each $P_{\xi}$ agrees with~$P$ on every cylinder. Hence,
each $%
\sigma$-additive $P_{\xi}$ must agree with $P$ on every Borel sets.
Because the sets $ \{D_{\xi}\dvtx \xi\in [ 0,1 ]  \} $ are
Borel and pairwise disjoint, at most one probability in $ \{P_{\xi
}\dvtx \xi\in [ 0,1 ]  \} $ agrees with $P$ on every Borel set.
Thus, there exists at most one $\sigma$-additive probability in $
\{P_{\xi}\dvtx \xi\in [ 0,1 ]  \} $. Therefore, the set
$ \{
P_{\xi}\dvtx \xi\in [ 0,1 ],P_{\xi}\neq P \} $, which is
included in $ \{ Q\in\Delta_{\mathrm{BD}}\dvtx Q_{\mathcal{F}}=P_{\mathcal
{F}}, Q\mbox{ is not }\sigma\mbox{-additive} \} $, has cardinality $%
\mathfrak{c}$. This completes the proof.
\end{pf}

\section{Strongly nonatomic probabilities}\label{sec9}

We now provide a technical result important for the proofs of Theorems~\ref{th2}
and~\ref{th5}. Throughout this subsection, $\Omega=\mathcal{X}^{\infty}$.

A $ \{ 0,1 \} $\emph{-probability} is a probability $Z\in
\mathbb{P}$ that satisfies $Z ( E ) \in \{ 0,1 \} $ for
every $%
E\in\Sigma$.

\begin{theorem}\label{th11}
Let $\mathcal{E}\subseteq\Sigma- \{ \varnothing \} $ be closed
under finite intersection. There exists a $ \{ 0,1 \} $-probability
$Z$ such that $Z ( E ) =1$ for every $E\in\mathcal{E}$.
\end{theorem}

\begin{pf}
This is a corollary of the ultrafilter theorem. See, for instance,
Aliprantis and Border (\citeyear{AliBor06}), Theorem 2.19.
\end{pf}

Every $P\in\mathbb{P}$ can be decomposed into a strongly nonatomic part
and a mixture of countably many $ \{ 0,1 \} $-probabilities.

\begin{theorem}[{[\citet{SobHam44}]}]\label{th12}
 For every $P\in\mathbb{P,}$ there exists an
opinion $P_{s}\in\Delta$ and a sequence $ ( Z_{i} )
_{i=1}^{\infty}$ of $ \{ 0,1 \} $-probabilities such that
\[
P=\alpha P_{s}+ ( 1-\alpha ) \sum_{i=1}^{\infty}
\beta_{i}Z_{i},
\]
where $\alpha,\beta_{i}\in [ 0,1 ] $ for every $i$ and $%
\sum_{i=1}^{\infty}\beta_{i}=1$.
\end{theorem}

Given an algebra $\mathcal{A}$, $P\in\mathbb{P} ( \mathcal{A}
) $
is \emph{strongly continuous} if for every $\varepsilon>0$ there exists a
partition $ \{ A_{1},\ldots,A_{n} \} $ of $\Omega$ such that $%
A_{i}\in\mathcal{A}$ and $P ( A_{i} ) <\varepsilon$ for every $i$.

\begin{theorem}\label{th13}
Let $\mathcal{A}$ be a $\sigma$-algebra. A probability $P\in\mathbb{P}
( \mathcal{A} ) $ is strongly continuous if and only if it is
strongly nonatomic.
\end{theorem}

\begin{pf}
See Bhaskara~Rao and Bhaskara~Rao (\citeyear{BhaBha83}), Theorem 11.4.5.
\end{pf}

\begin{theorem}\label{th14}
For every $P\in\mathbb{P,}$ there exists an opinion $\widetilde{P}\in
\Delta
$ such that
\[
P ( U ) \leq\widetilde{P} ( U )
\]
for every open set $U$.
\end{theorem}

\begin{pf}
We first prove the result for the case where $P$ is a $ \{ 0,1
\} $%
-probability. To this end, fix a $ \{ 0,1 \} $-probability
$Z$. Let
$\mathcal{C=} \{ U\dvtx U\mbox{ open},\break Z ( U ) =1 \} $. The
collection $\mathcal{C}$ is closed under finite intersection.

We now construct a sequence $ ( D_{n} ) _{n=1}^{\infty}$ of
countable, dense and pairwise disjoint subsets of $\Omega$. The proof of
this claim proceeds by induction. The space $\Omega=\mathcal
{X}^{\infty}$
is separable, so it has a countable dense subset $D_{1}$. Assume that for
some $N$ the sets $D_{1},\ldots,D_{N}$ have been defined and satisfy the
desired properties. Let $ ( V_{k} ) _{k=1}^{\infty}$ be a
countable base of $\Omega$ and for each $k$ pick a path $\omega_{k}\in
V_{k}-\bigcup_{n=1}^{N}D_{n}$. Let $D_{N+1}= \{ \omega_{1},\omega
_{2},\ldots \} $. This completes the induction step and the proof of the
claim.

For every $n$, the collection
\[
\mathcal{C}_{n}= \{ U\cap D_{n}\dvtx U\in\mathcal{C} \}\vadjust{\goodbreak}
\]
is closed under finite intersection and does not contain the empty set. By
Theorem~\ref{th11}, for every $n$ there exists a $ \{ 0,1 \}
$-probability $%
Z_{n}$ such that $Z_{n} ( E ) =1$ for every $E\in\mathcal{C}_{n}$.
For each $n$ and for every open set $U$, if $Z ( U ) =1$ then $
U\cap D_{n}\in\mathcal{C}_{n}$ and $Z_{n} ( U ) =1$. Hence, $%
Z ( U ) \leq Z_{n} ( U ) $. Let $\lambda$ be a strongly
continuous finitely additive probability on $ (
\mathbb{N},2^{\mathbb{N}
} ) $ and define the function~$\widetilde{Z}\dvtx \Sigma\rightarrow
[
0,1 ] $ as
\[
\widetilde{Z} ( E ) =\int_{%
\mathbb{N}
}Z_{n} ( E ) \,d
\lambda ( n )
\]
for every $E\in\Sigma$.

It follows from the additivity of the integral that $\widetilde{Z}\in
\mathbb{P}$. For every open set $U$, if $Z ( U ) =1$ then $%
\widetilde{Z} ( U ) =\int_{%
\mathbb{N}
}Z_{n} ( U ) \,d\lambda ( n ) =\int_{%
\mathbb{N}
}1\,d\lambda ( n ) =1$.\vspace*{1pt} Hence, $Z ( U ) \leq\widetilde
{Z}%
( U ) $. It remains to prove that $\widetilde{Z}$ is strongly
nonatomic. By Theorem \ref{th13}, it is enough to prove it is strongly continuous.
Fix $\varepsilon>0$. Since $\lambda$ is strongly continuous, we can find a
partition $ \{ \Pi_{1},\ldots,\Pi_{k} \} $ of $%
\mathbb{N}
$ such that $\lambda ( \Pi_{i} ) \leq\varepsilon$ for every $%
i=1,\ldots,k$. Consider now the partition $ \{ \bigcup_{m\in\Pi
_{1}}D_{m},\ldots,\bigcup_{m\in\Pi_{k}}D_{m},\vspace*{1pt}\break  ( \bigcup_{m\in
\mathbb{N}
}D_{m} ) ^{c} \} $ of $\Omega$. For every $i=1,\ldots,k$ and
$n$, we
have $Z_{n} ( \bigcup_{m\in\Pi_{i}}D_{m} ) =1_{\Pi_{i}} (
n ) $ where $1_{\Pi_{i}}$ is the indicator function of $\Pi_{i}$.
Therefore,
\begin{eqnarray*}
\widetilde{Z} \biggl( \bigcup_{m\in\Pi_{i}}D_{m}
\biggr) &=&\int_{%
\mathbb{N}
}Z_{n} \biggl( \bigcup
_{m\in\Pi_{i}}D_{m} \biggr) \,d \lambda ( n )
\\
&=&\int_{%
\mathbb{N}
}1_{\Pi_{i}} ( n ) \,d\lambda ( n )
\\
&=&\lambda ( \Pi_{i} ) \leq\varepsilon
\end{eqnarray*}
and $\widetilde{Z} (  ( \bigcup_{m\in
\mathbb{N}
}D_{m} ) ^{c} ) =0$. This proves that $\widetilde{Z}$ is strongly
continuous.

Now let $P$ be any finitely additive probability. By the Hammer--Sobczyk
decomposition, we can write $P$ as the convex combination
\[
P=\alpha P_{s}+ ( 1-\alpha ) \sum_{i=1}^{\infty}
\beta_{i}Z_{i},
\]
where $P_{s}$ is strongly nonatomic and each $Z_{i}$ is a $ \{
0,1 \} $-probability. For each $i$, let $\widetilde{Z}_{i}$ be an
opinion such that $Z_{i} ( U ) \leq\widetilde{Z}_{i} (
U ) $ for  every open set $U$. Define $\widetilde{P}=\alpha
P_{s}+ ( 1-\alpha ) \sum_{i=1}^{\infty}\beta\widetilde{Z}_{i}$.
It is easy to see that $\widetilde{P}$ is strongly continuous. By Theorem
\ref{th13}, it is an opinion. For every open set $U$, we have
\begin{eqnarray*}
P ( U ) &=&\alpha P_{s} ( U ) + ( 1-\alpha ) \sum
_{i=1}^{\infty}\beta_{i}Z_{i} ( U )
\\
&\leq&\alpha P_{s} ( U ) + ( 1-\alpha ) \sum
_{i=1}^{\infty}\beta_{i}\widetilde{Z}_{i}
( U ) =\widetilde{P}%
( U )
\end{eqnarray*}
as desired.\vadjust{\goodbreak}
\end{pf}

\section{Proofs of Theorems \texorpdfstring{\lowercase{\protect\ref{th2}}--\lowercase{\protect\ref{th6}}}{2--6}}\label{sec10}

\begin{theorem}[{[\citet{Fan53}]}]\label{th15}
Let $X$ and $Y$ be convex subsets of two vector spaces. Let $f\dvtx X\times
Y\rightarrow
\mathbb{R}
$. If $X$ is compact Hausdorff and $f$ is concave with respect to $Y$ and
convex and lower semicontinuous with respect to $X$, then
\[
\sup_{y\in Y}\min_{x\in X}f ( x,y ) =\min
_{x\in X}\sup_{y\in
Y}f ( x,y ).
\]
\end{theorem}

See \citet{Fan53} for a more general version of this theorem.\vspace*{-2pt}

\begin{pf*}{Proof of Theorem \ref{th2}}
Define the function $V\dvtx \mathbb{P}\times\Delta_{f}\Delta
\longrightarrow
\mathbb{R}
$ as
\[
V ( P,\zeta ) =\int\zeta \bigl( \bigl\{ Q\in\Delta\dvtx \omega \notin T ( Q )
\bigr\} \bigr) \,dP ( \omega )
\]
for every $ ( P,\zeta ) \in\mathbb{P}\times\Delta
_{f}\Delta$.
The function $V$ is affine in each variable and continuous with respect
to $%
\mathbb{P}$. The weak* topology is Hausdorff. Moreover, it follows from the
Riesz representation and Banach--Alaoglu theorems that $\mathbb{P}$ is
compact [see Aliprantis and Border (\citeyear{AliBor06}),
Theorems 14.4 and 6.21]. All the
conditions of Fan's Minmax theorem are verified, therefore,
%
\begin{equation}
\sup_{\zeta\in\Delta_{f}\Delta}\min_{P\in\mathbb{P}}V ( P,\zeta ) =\min
_{P\in\mathbb{P}}\sup_{\zeta\in\Delta_{f}\Delta}V ( P,\zeta ).
\end{equation}

By Theorem \ref{th14}, for every $P\in\mathbb{P}$ there exists $\widetilde
{P}\in
\Delta$ such that $P ( U ) \leq\widetilde{P} ( U )
$ for
every open set $U$. Because $T(\widetilde{P})$ is an open set, then
$P(T(%
\widetilde{P}))\leq\widetilde{P}(T(\widetilde{P}))\leq\varepsilon$.
That is,
$V(P,\delta_{\widetilde{P}})=1-P(T(\widetilde{P}))\geq1-\varepsilon$. Thus,
\[
\sup_{\zeta\in\Delta_{f}\Delta}\min_{P\in\mathbb{P}}V ( P,\zeta ) =\min
_{P\in\mathbb{P}}\sup_{\zeta\in\Delta_{f}\Delta}V ( P,\zeta ) \geq\min
_{P\in\mathbb{P}}V(P,\delta_{\widetilde
{P}})\geq 1-\varepsilon.
\]
For every $\delta\in(0,1-\varepsilon]$, there exists a strategy $\zeta
$%
such that $V ( P,\zeta ) >1-\varepsilon-\delta$ for every $P\in
\mathbb{P}$. In particular,
\[
V ( \delta_{\omega},\zeta ) =\zeta \bigl( \bigl\{ Q\in \Delta \dvtx \omega
\notin T ( Q ) \bigr\} \bigr) \geq1-\varepsilon-\delta
\]
for every path $\omega$.\vspace*{-2pt}
\end{pf*}

\begin{pf*}{Proof of Theorem \ref{th3}}
Fix $\varepsilon>0$ and a path $\omega$. By Theorem \ref{th9}, for every $P\in
\Delta
_{\mathrm{BD}}$ we can choose a partition $ \{ C_{1},\ldots,C_{n} \} $ of
cylinders such that $P ( C_{i} ) <\varepsilon$ for every $i=1,\ldots,n$.
Let $\omega\in C_{i}$. There exists a time $t_{P}$ such that $\omega
^{t_{P}}=C_{i}$. Hence, $P ( \omega^{t_{P}} ) <\varepsilon$.
Define a
test $T$ as $T ( P ) =\omega^{t_{P}}$ for every opinion $P\in
\Delta_{\mathrm{BD}}$. The test $\Delta_{\mathrm{BD}}$-controls for type I errors with
probability $1-\varepsilon$.\looseness=-1

Now, let $ ( P_{1},\ldots,P_{n} ) $ be the support of a strategy $
\zeta$. Choose a time $t$ such that $t\geq t_{P_{i}}$ for each $i=1,\ldots,n$.
Then
\[
\omega^{t}\subseteq\bigcap_{i=1}^{n}T
( P_{i} )
\]
hence, $\zeta (  \{ P\in\Delta\dvtx \tilde{\omega}\notin T (
P )  \}  ) =0$ for every path $\tilde{\omega}$ in $\omega
^{t} $.\vadjust{\goodbreak}
\end{pf*}

\begin{pf*}{Proof of Theorem \ref{th4}}
Recall that $\mathcal{F}$ is the algebra generated by cylinders. Fix a
probability $\pi\in\mathbb{P} ( \mathcal{F} ) $ such that
$\pi
( \omega^{t} ) \rightarrow0$ as $t\rightarrow\infty$ for every
path $\omega$. Now let $\mathcal{I}$ be a strictly proper ideal and
consider the collection of events
%
\begin{equation}\label{eqD2}
\mathcal{A}= \bigl\{ ( F\cap L ) \cup S\dvtx F\in\mathcal{F}, S\in\mathcal{I}
\cap \Sigma, L^{c}\in\mathcal{I}\cap\Sigma \bigr\}.
\end{equation}
We prove it is an algebra. If $ ( F\cap L ) \cup S\in
\mathcal{A,}$ then its complement is equal to
\[
( F\cap L ) ^{c}\cap S^{c}= \bigl( F^{c}\cup
L^{c} \bigr) \cap S^{c}= \bigl( F^{c}\cap
S^{c} \bigr) \cup \bigl( L^{c}\cap S^{c} \bigr)
\]
since $S^{c}$ is large and $L^{c}\cap S^{c}$ is small we have that
$ (
F^{c}\cap S^{c} ) \cup ( L^{c}\cap S^{c} ) \in\mathcal{A}$.
Using the notation in (\ref{eqD2}), let $ ( F_{1}\cap L_{1} ) \cup S_{1}$
and $ ( F_{2}\cap L_{2} ) \cup S_{2}$ belong to $\mathcal{A}$.
Observe that $L=L_{1}\cap L_{2}$ is large and, therefore $L_{1}-L$ and $
L_{2}-L$ are small. We can write
%
\begin{eqnarray}\label{eqD3}
( F_{1}\cap L_{1} ) \cup S_{1}\cup (
F_{2}\cap L_{2} ) \cup S_{2} &=& (
F_{1}\cap L ) \cup ( F_{2}\cap L ) \cup S
\nonumber
\\[-8pt]
\\[-8pt]
\nonumber
&=& \bigl( ( F_{1}\cup F_{2} ) \cap L \bigr) \cup S,
\nonumber
\end{eqnarray}
where $S= ( F_{1}\cap ( L_{1}-L )  ) \cup S_{1}\cup
( F_{2}\cap ( L_{2}-L )  ) \cup S_{2}$ is a union of
small sets. This proves that $ ( F_{1}\cap L_{1} ) \cup
S_{1}\cup
( F_{2}\cap L_{2} ) \cup S_{2}\in\mathcal{A}$. We conclude
that $%
\mathcal{A}$ is an algebra. By construction, $\mathcal{F}\subseteq
\mathcal{A%
}\subseteq\Sigma$.

Define a set function $\widetilde{\pi}\dvtx \mathcal{A}\rightarrow [
0,1%
] $ as
\[
\widetilde{\pi} \bigl( ( F\cap L ) \cup S \bigr) =\pi ( F )
\]
for each $ ( F\cap L ) \cup S\in\mathcal{A}$.

We verify that $\widetilde{\pi}$ is well defined. Using the notation
in (\ref{eqD2}),
let $ ( F_{1}\cap L_{1} ) \cup S_{1}= ( F_{2}\cap
L_{2} )
\cup S_{2}$. Equivalently,
\[
( F_{1}\cap F_{2}\cap L_{1} ) \cup \bigl(
F_{1}\cap F_{2}^{c}\cap L_{1} \bigr)
\cup S_{1}= ( F_{2}\cap L_{2} ) \cup
S_{2}.
\]
Therefore, $F_{1}\cap F_{2}^{c}\cap L_{1}\subseteq S_{2}$. Hence,
$F_{1}\cap
F_{2}^{c}\cap L_{1}$ is small. But also $F_{1}\cap F_{2}^{c}\cap
L_{1}^{c}\subseteq L_{1}^{c}$ is small, hence $F_{1}\cap F_{2}^{c}$ is
small. The set $F_{1}\cap F_{2}^{c}$ is either empty or a union of
cylinders. By the definition of strictly proper ideal $F_{1}\cap F_{2}^{c}$
must be empty. Similarly, $F_{2}\cap F_{1}^{c}=\varnothing$. Hence, $%
F_{1}=F_{2}$, and $\widetilde{\pi} (  ( F_{1}\cap L_{1} )
\cup
S_{1} ) =\widetilde{\pi} (  ( F_{2}\cap L_{2} ) \cup
S_{2} ) $.

We prove $\widetilde{\pi}$ is additive. Let $ ( F_{1}\cap
L_{1} )
\cup S_{1}$ and $ ( F_{2}\cap L_{2} ) \cup S_{2}$ be two disjoint
sets belonging to $\mathcal{A}$. The sets $F_{1}$ and $F_{2}$ are disjoint.
To see this, notice that $F_{1}\cap F_{2}\cap L_{1}\cap L_{2}=\varnothing$
implies $F_{1}\cap F_{2}\subseteq ( L_{1}\cup L_{2} ) ^{c}$. The
set $F_{1}\cap F_{2}$ is either empty or a union of cylinders. Since
$ (
L_{1}\cup L_{2} ) ^{c}$ is small, it must be empty. Let
$L=L_{1}\cap
L_{2}$  and $S= ( F_{1}\cap ( L_{1}-L )  ) \cup
S_{1}\cup ( F_{2}\cap ( L_{2}-L )  ) \cup S_{2}$.  Similar to (\ref{eqD3}), we have
\begin{eqnarray*}
\widetilde{\pi} \bigl( ( F_{1}\cap L_{1} ) \cup
S_{1}\cup ( F_{2}\cap L_{2} ) \cup
S_{2} \bigr) &=&\widetilde{\pi} \bigl( \bigl( ( F_{1}\cup
F_{2} ) \cap L \bigr) \cup S \bigr)
\\
&=&\pi ( F_{1}\cup F_{2} )
\\
&=&\pi ( F_{1} ) +\pi ( F_{2} )
\\
&=&\widetilde{\pi} \bigl( ( F_{1}\cap L_{1} ) \cup
S_{1} \bigr) +%
\widetilde{\pi} \bigl( ( F_{2}
\cap L_{2} ) \cup S_{2} \bigr).
\end{eqnarray*}
Therefore, $\widetilde{\pi}$ is a finitely additive probability
defined on $%
( \Omega,\mathcal{A} ) $. By construction, it satisfies $%
\widetilde{\pi} ( S ) =0$ for every $S\in\mathcal{I}\cap
\Sigma$.

Consider the set of extensions $E ( \widetilde{\pi},\mathcal
{A},\Sigma
) $ and let $P$ be one of its extreme points. We prove that $P$
is an
extreme point of $E ( \pi,\mathcal{F},\Sigma ) $. Write $P$
as $%
P=\alpha Q+ ( 1-\alpha ) R$ with $Q,R\in E ( \pi,\mathcal
{F}%
,\Sigma ) $. Let $\widetilde{\pi}_{Q}$ and $\widetilde{\pi
}_{R}$ be
the restriction of $Q$ and $R$ on $\mathcal{A}$. Since $P$ is an extension
of $\widetilde{\pi}$, we have $\widetilde{\pi}=\alpha\widetilde{\pi}
_{Q}+ ( 1-\alpha ) \widetilde{\pi}_{R}$. We claim that $%
\widetilde{\pi}=\widetilde{\pi}_{Q}=\widetilde{\pi}_{R}$. For every
$S\in
\mathcal{I}\cap\Sigma$, since $\widetilde{\pi} ( S ) =0$,
then $%
\widetilde{\pi}_{Q} ( S ) =\widetilde{\pi}_{R} ( S
) =0$%
. Therefore,
\[
\widetilde{\pi}_{Q} \bigl( ( F\cap L ) \cup S \bigr) =
\widetilde{%
\pi}_{Q} ( F\cap L ) =\widetilde{
\pi}_{Q} ( F ) =\pi ( F ) =\widetilde{\pi} \bigl( ( F\cap L ) \cup S
\bigr)
\]
for every event $  ( F\cap L ) \cup S\in\mathcal{A}$. The
same is
true for $\widetilde{\pi}_{R}$. Therefore, $\widetilde{\pi
}_{Q}=\widetilde{%
\pi}_{R}=\widetilde{\pi}$. This proves that $Q,R\in E ( \widetilde
{\pi
},\mathcal{A},\Sigma ) $. Because $P$ is an extreme point of
$E (
\widetilde{\pi},\mathcal{A},\Sigma ) $, then $P=Q=R$. This concludes
the proof that $P$ is an extreme point of $E ( \pi,\mathcal
{F},\Sigma
) $. By Theorem \ref{th8}, $P$ satisfies the Blackwell--Dubins property.

It remains to prove that $P$ is strongly nonatomic. Since $\pi (
\omega^{t} ) \rightarrow0$ for every $\omega$, $\pi$ is strongly
continuous. A fortiori, $P$ is strongly continuous and also strongly
nonatomic by Theorem \ref{th13}.
\end{pf*}

\begin{pf*}{Proof of Theorem \ref{th5}}
Define the function $V\dvtx \mathbb{P}\times\Delta_{f}\Delta^{\ast
}\rightarrow\mathbb{R}$ as
\[
V ( P,\zeta ) =\int\zeta \bigl( \bigl\{ Q\in\Delta^{\ast
}\dvtx \omega
\notin T ( Q ) \bigr\} \bigr) \,dP ( \omega )
\]
for all $ ( P,\zeta ) \in\mathbb{P}\times\Delta_{f}\Delta
^{\ast}$. Given $P\in\mathbb{P}$, by Theorem \ref{th14} there exists an
opinion $Q$
such that $P ( U ) \leq Q ( U ) $ for every open set $U$.
By Theorem 4 in \citet{Reg85}, we can find a conditional opinion
$Q^{\ast
}\in\Delta^{\ast}$ such that $Q=Q^{\ast} ( \cdot|\Omega
) $.
Then $P ( T ( Q^{\ast} )  ) \leq Q ( T (
Q^{\ast
} )  ) =Q^{\ast} ( T ( Q^{\ast} )  ) \leq
\varepsilon$. The proof is complete by replicating the argument used in the
proof of Theorem \ref{th2}.
\end{pf*}

\begin{pf*}{Proof of Theorem \ref{th6}}
We first prove that for every $P\in\Delta_{\mathrm{BD}}^{\ast}$ and every
path $%
\omega$, $\lim_{t}P ( \omega^{t} ) =0$. We argue by
contradiction. Let $\omega_{o}$ be a path such that $\inf_{t}P (
\omega
_{o}^{t} ) =\delta>0$. Fix a sequence of positive real numbers
$ (
\xi_{t} ) $ such that
\[
P \bigl( \omega_{o}^{t} \bigr) =\delta+\xi_{t}
\]
for every $t$.

Fix $\varepsilon\in ( 0,\frac{1}{2} ) $. Because $P$ is strongly
nonatomic, for every time $t$ we can find an event $F^{t}\subseteq$ $%
\omega_{0}^{t}$ such that $P ( F^{t} ) =\frac{1}{2}P (
\omega_{0}^{t} ) $. For every $n$, we have
\begin{eqnarray*}
P \bigl( F^{t}|\omega_{0}^{t+n} \bigr) &=&
\frac{P ( F^t \cap\omega
_{o}^{t+n} ) }{P ( \omega_{o}^{t+n} ) }
\\
&= &\frac{P ( F^{t}  ) -P ( \omega
_{o}^{t}-\omega_{o}^{t+n} ) }{\delta+\xi_{t+n}}
\\
&=&\frac{({1}/{2}) ( \delta+\xi_{t} ) - ( \xi_{t%
}-\xi_{t+n} ) }{\delta+\xi_{t+n}}=\frac{1}{2}\frac{\delta+\xi
_{t%
}}{\delta+\xi_{t+n}}-\frac{\xi_{t}-\xi_{t+n}}{\delta+\xi_{t+n}}.
\end{eqnarray*}
We can therefore fix $\bar{t}$ large enough such that $F=F^{\bar{t}}$
satisfies $P ( F|\omega_{0}^{t} )
\in ( \frac{1}{2}-\varepsilon,\frac{1}{2}+\varepsilon ) $ for
every $t>%
\bar{t}$.

Let $Q$ be the opinion defined as
\[
Q ( E ) =\frac{P ( E\cap F ) }{P ( F ) }
\]
for every event $E$. Then $Q\ll P ( \cdot|\Omega ) $. By Theorem
4 in \citet{Reg85}, we can find a conditional opinion $\widetilde{Q}$
satisfying $\widetilde{Q} ( \cdot|\Omega ) =Q$. The proof of
the claim will be
concluded by showing that $P$ does not merge with $\widetilde{Q}$. Note
that for every $t>\bar{t}$
\[
\widetilde{Q} \bigl( \omega_{0}^{t} \bigr) =
\frac{P ( F\cap\omega
_{0}^{t} ) }{P ( \omega_{0}^{t} ) }\frac{P ( \omega
_{0}^{t} ) }{P ( F ) }=P \bigl( F|\omega_{0}^{t}
\bigr) \frac{%
P ( \omega_{0}^{t} ) }{({1}/{2})P ( \omega_{0}^{\bar{t}
} ) }=P \bigl( F|\omega_{0}^{t} \bigr) 2
\frac{\delta+\xi
_{t}}{\delta
+\xi_{\bar{t}}}
\]
hence, for all $t>\bar{t}$
\[
\widetilde{Q} \bigl( \omega_{0}^{t} \bigr) \geq ( 1-2
\varepsilon ) \frac{\delta}{\delta+\xi_{\bar{t}}}.
\]
Moreover, for every $t>\bar{t,}$
\begin{eqnarray*}
\widetilde{Q} \biggl( \biggl\{ \omega\dvtx \sup_{E} \bigl
\llvert Q \bigl( E|\omega ^{t} \bigr) -P \bigl( E|\omega^{t}
\bigr) \bigr\rrvert >\frac
{1}{2}-\varepsilon \biggr\} \biggr) &>&\widetilde{Q}
\biggl( \biggl\{ \omega\dvtx P \bigl( F|\omega ^{t} \bigr) <
\frac{1}{2}+\varepsilon \biggr\} \biggr)
\\
&\geq&\widetilde{Q} \bigl( \omega_{0}^{t} \bigr),
\end{eqnarray*}
where the first equality follows from $\widetilde{Q} ( F|\omega
_{0}^{t} ) =1$ and the second equality follows from $P (
F|\omega
_{0}^{t} ) <\frac{1}{2}+\varepsilon$. Because the sequence $ (
\widetilde{Q} ( \omega_{0}^{\bar{t}} ),\widetilde{Q} (
\omega
_{0}^{\bar{t}+1} ),\ldots ) $ is bounded away from $0$, $P$
does not
merge to $\widetilde{Q}$.
Therefore, we can conclude that for every $P\in\Delta_{\mathrm{BD}}^{\ast}$
and every path $%
\omega$, $\lim_{t}P ( \omega^{t} ) =0$.

Now fix a path $\omega$ and $\varepsilon>0$. We can find for every $P\in
\Delta_{\mathrm{BD}}^{\ast}$ a time $t_{P}$ such that $P ( \omega
^{t_{P}} ) <\varepsilon$. Because $\mathcal{X}$ is endowed with the
discrete topology, $\omega^{t_{P}}$ is an open set. Let $T (
P )
=\omega^{t_{P}}$ for every $P\in\Delta_{\mathrm{BD}}^{\ast}$. The test
$\Delta
_{\mathrm{BD}}^{\ast}$-controls for type~I error with probability $1-\varepsilon
$. The
same argument in the proof of Theorem \ref{th3} shows that $T$ is nonmanipulable.
\end{pf*}
\end{appendix}
\section*{Acknowledgements}
We thank Ehud Kalai, Wojciech Olszewski, Eran Shmaya, Marciano Siniscalchi
and Rakesh Vohra for useful discussions. We are grateful to the Editor and
the referees for their thoughtful comments, for simplifying Example \ref{ex1} and
for stimulating the results in Section~\ref{sec6}. We also thank the seminar
audiences at the Fifth Transatlantic Theory Workshop, the Summer
meeting of
the Econometric Society 2012, XIII Latin American Workshop in Economic
Theory, the 4th Workshop on Stochastic Methods in Game Theory, the
Washington University seminar series and the Paris Game Theory Seminar.
All errors are ours.

%



\printaddresses

\end{document}